\documentclass{article}

\usepackage{mib}

\usepackage{color}
\usepackage[utf8]{inputenc} 
\usepackage[T1]{fontenc}    
\usepackage{hyperref}       
\usepackage{url}            
\usepackage{booktabs}       
\usepackage{amsfonts}       
\usepackage{nicefrac}       
\usepackage{microtype}      
\usepackage{lipsum}
\usepackage{graphicx}
\usepackage{color}
\usepackage{tikz}
\usepackage{amsthm,amsmath,amssymb}
\usepackage[english, russian]{babel}
\title{Description of fixed points of an infinite dimensional operator}

\author{ Olimov Umrbek \\
  V.I. Romanovskiy Institute of Mathematics\\
   Uzbekistan Academy of Sciences\\
    Tashkent, Uzbekistan\\
     umrbek.olimov.92@mail.ru\\
}

\begin{document}
\maketitle
\begin{center}
	
	{\large \sc \bf Abstract}
\end{center}
\begin{abstract}
{ We consider an infinite-dimensional non-linear operator related to a hard core (HC) model with a countable set  $\mathbb{N}$ of spin values. It is known that finding the fixed points of an infinite-dimensional operator is generally impossible. But we have fully analyzed the fixed points of an infinite-dimensional operator by applying a technique of reducing an  infinite-dimensional operator to a two-dimensional operator. The set of parameters is divided into subsets $A_{i,j},$ where the index $i$ means the number of fixed points on the line $y=x$, $j$ means the number of fixed points outside of $y=x.$ The number of fixed points can be up to seven, and the explicit form of each fixed point is found.}

\end{abstract}
\label{firstpage}
\bigskip

\label{firstpage}

\lhead{{\sf Olimov U.~R.~ Description of fixed points of an operator}}

{\bf \large 1. Introduction}

Consider an operator $F$ defined on a topological space $X$.

If $F(x)=x$ then the point $x\in X$ is called a fixed point for $F$.

In this paper we are interested to find fixed points of an infinite-dimensional non-linear operator mapping strictly positive sequences from $l_1$ to itself.

 Denote
$$l^{1}_{+}=\left\{x=(x_{1},x_{2},\dots,x_{n},\dots) \, : \, x_{i}>0, \, \|x\|=\sum_{j=1}^{\infty}x_j<\infty\right\}.$$

In this paper we study fixed points of the operator $F$ defined on $ l^{1}_{+}$ as follows
\begin{equation}\label{e4o}
F:	\begin{array}{lll}
		x'_{2n-1}=\lambda_{2n-1}\cdot \left(\dfrac{1+\sum_{j=1}^{\infty}x_{2j-1}}{1+\theta+\|x\|}\right)^2,\\[2mm]
		x'_{2n}=\lambda_{2n}\cdot \left(\dfrac{1+\sum_{j=1}^{\infty}x_{2j}}{1+\theta+\|x\|}\right)^2,
	\end{array}
\end{equation}
where $n=1,2,\dots$, $\theta>0$ and $\lambda=(\lambda_1, \lambda_2, \dots)\in l^{1}_{+}$.

Note that these fixed points are very important to study Gibbs measures of a HC-model (see \cite{OR}).

The paper is organized as follows. In Section 2 we reduce our operator to a two-dimensional operator and show that finding fixed points of operator (\ref{e4o}) is equivalent to finding of such points for the two-dimensional operator.  In Section 3  we find some invariant sets and give full description of fixed points, the number of which is up to seven. Namely, we divide the set of parameters to seven subsets, for parameters from each set the number of fixed points is given (up to seven).

{\bf \large 2. Restriction to two-dimensional operator}

In this section we reduce our problem to a two-dimensional case.

\textbf{Lemma 1.} If $\lambda=(\lambda_1, \lambda_2, \dots)\in l^{1}_{+}$ then $F$ maps $l^{1}_{+}$ to itself.

{\sf Proof. }
Follows from the straightforward inequalities: $x'_i<\lambda_i$	for any $i=1,2,\dots.$

Recall definition of dynamical system \cite{Romb}:

Let $F:X\rightarrow X$ and $\forall x^{(0)}\in X$ be given. Then the sequence $x^{(0)},x^{(1)}=F(x^{(0)}),x^{(2)}=F(x^{(1)}),...$ is called a dynamical system (DS).

Define two-dimensional operator
 $W: z=(x,y)\in \mathbb R^2_+\to z'=(x',y')=W(z)\in \mathbb R^2_+$ by
\begin{equation}\label{op}
	W: \ \ \begin{array}{ll}
		x'=L_1\cdot \left(\dfrac{1+x}{1+\theta+x+y}\right)^2\\[2mm]
		y'=L_2\cdot\left(\dfrac{1+y}{1+\theta+x+y}\right)^2,
	\end{array}
\end{equation}
where $\theta>0$ and $L_i>0$ are parameters.

\textbf{Lemma 2.} The infinite dimensional DS generated by the operator $F$ is fully represented by the two-dimensional DS generated by the operator $W$.

{\sf Proof. }

Introduce
$$M_1=\sum_{j=1}^{\infty}x_{2j-1}, \ \ M_2=\sum_{j=1}^{\infty}x_{2j},$$ $$L_1=\sum_{j=1}^{\infty}\lambda_{2j-1}, \ \ L_2=\sum_{j=1}^{\infty}\lambda_{2j}.$$
Note that $\|x\|=M_1+M_2$.

Then the operator (\ref{e4o}) is reduced to
\begin{equation}\label{e4}
	\begin{array}{lll}
	x'_{2n-1}=\lambda_{2n-1}\left(\dfrac{1+M_1}{1+\theta+\|x\|}\right)^2,\\[2mm]
		x_{2n}'=\lambda_{2n}\left(\dfrac{1+M_2}{1+\theta+\|x\|}\right)^2.
	\end{array}
\end{equation}
Denoting  $$M_1'=\sum_{j=1}^{\infty}x'_{2j-1}, \ \ M_2'=\sum_{j=1}^{\infty}x'_{2j},$$
and summing the equalities (\ref{e4}) we get
\begin{equation}\label{ef4}
	\begin{array}{ll}
		M'_1=L_1\left(\dfrac{1+M_1}{1+\theta+M_1+M_2}\right)^2\\[2mm]
		M'_2=L_2\left(\dfrac{1+M_2}{1+\theta+M_1+M_2}\right)^2.
	\end{array}
\end{equation}
This is the operator (\ref{op}).
Now it is easy to see that each fixed point $(M_1, M_2)$ to (\ref{ef4}), by formula (\ref{e4}), uniquely defines a fixed point of operator (\ref{e4o}).

Moreover, let $x^{(m)}=F^m(x^{(0)})$, $m\geq 1$ be trajectory of the arbitrary point $x^{(0)}\in l^1_+$, i.e., coordinates of which satisfy
\begin{equation}\label{etr}
	\begin{array}{lll}
		x^{(m+1)}_{2n-1}=\lambda_{2n-1}\cdot \left(\dfrac{1+\sum_{j=1}^{\infty}x^{(m)}_{2j-1}}{1+\theta+\|x^{(m)}\|}\right)^2=\lambda_{2n-1}\cdot \left(\dfrac{1+M_1^{(m)}}{1+\theta+M_1^{(m)}+M_2^{(m)}}\right)^2,\\[2mm]
		x^{(m+1)}_{2n}=\lambda_{2n}\cdot \left(\dfrac{1+\sum_{j=1}^{\infty}x^{(m)}_{2j}}{1+\theta+\|x^{(m)}\|}\right)^2=\lambda_{2n}\cdot \left(\dfrac{1+M_2^{(m)}}{1+\theta+M_1^{(m)}+M_2^{(m)}}\right)^2,
	\end{array}
\end{equation}

Summing the last equalities we get
\begin{equation}\label{me}
	\begin{array}{ll}
		M^{(m+1)}_1=L_1\left(\dfrac{1+M^{(m)}_1}{1+\theta+M^{(m)}_1+M^{(m)}_2}\right)^2\\[2mm]
		M^{(m+1)}_2=L_2\left(\dfrac{1+M^{(m)}_2}{1+\theta+M^{(m)}_1+M^{(m)}_2}\right)^2.
	\end{array}
\end{equation}
Using now (\ref{me}) rewrite (\ref{etr}) as
\begin{equation}\label{er}
	\begin{array}{lll}
		x^{(m+1)}_{2n-1}=\lambda_{2n-1}\cdot {M_1^{(m+1)}\over L_1},\\[2mm]
		x^{(m+1)}_{2n}=\lambda_{2n}\cdot {M_2^{(m+1)}\over L_2}.
	\end{array}
\end{equation}
Thus if we know behavior of $(M^{(m)}_1, M^{(m)}_2)$, $m=0,1,\dots$ given by (\ref{me}) then by formula (\ref{er}) we can describe behavior of $x^{(m)}$.   For example, if $\lim_{m\to \infty}(M^{(m)}_1, M^{(m)}_2)=(a,b)$ then
\begin{equation}\label{xab}\lim_{m\to \infty}x^{(m)}=\left({a\over L_1}\lambda_1, {b\over L_2}\lambda_2, {a\over L_1}\lambda_3, {b\over L_2}\lambda_4, \dots\right).
\end{equation} 	
Lemma is proved.

\textbf{Remark 1.}By Lemma 2 we reduced investigation of trajectories $(x^{(m)})_{m=0}^\infty$, for all $x^{(0)}\in l^1_+$ to the study trajectories generated by the two-dimensional non-linear operator $W$ given in (\ref{op}). In the sequel we study fixed points of $W$. For simplicity we assume $L_1=L_2=L$.

{\bf \large 3. Fixed points for the case $L_1=L_2=L$}

Consider operator (\ref{op}) for the case $L_1=L_2=L$.

\textbf{3.1. Invariant sets.}
Denote
$$M_{-}=\{(x,y)\in \mathbb R^2_+: x<y\},$$
$$M_0=\{(x,y)\in \mathbb R^2_+: x=y\},$$
$$M_+=\{(x,y)\in \mathbb R^2_+: x>y\}.$$
\textbf{Lemma 3.}  If $L_1=L_2=L$ then the sets $M_\epsilon$, $\epsilon=-,0,+$ are invariant with respect to operator $W$, i.e., $W(M_\epsilon)\subset M_\epsilon$.

{\sf Proof. } Follows from the following equality
$$x'-y'=L \dfrac{(1+x)^2-(1+y)^2}{(1+\theta+x+y)^2}=
L \dfrac{(x-y)(2+x+y)}{(1+\theta+x+y)^2}.$$

\textbf{3.2. On the invariant set $M_0$.}
Reduce the operator $W$ defined by (\ref{op}) on the invariant set $M_0$, then we get
\begin{equation}\label{f}
x'=f(x):=L\left(\frac{1+x}{1+\theta+2x}\right)^2.
\end{equation}
Denote
$$\hat L_1=\frac{2\theta^2+76\theta-142-(2\theta-34)\sqrt{\theta^2-18\theta+17}}{16},$$

$$\hat L_2=\frac{2\theta^2+76\theta-142+(2\theta-34)\sqrt{\theta^2-18\theta+17}}{16}.$$
In \cite[Section 4.1]{OR} it is shown that
\begin{itemize}
\item[1)] If $\theta\in(0, 17]$, $L>0$ or $\theta>17$, $L\notin(\hat L_1, \hat L_2)$, then $f$ (defined in (\ref{f})) has unique fixed point, denoted by $x_1$;

\item[2)] If $\theta>17$ and $L=\hat L_1$ or $L=\hat L_2$ then $f$ has two fixed points, denoted by $x_1$, $x_2$ with $x_1<x_2$;

\item[3)] If $\theta>17$, $L\in(\hat L_1, \hat L_2)$ then $f$ has three fixed points $x_i$, $i=1, 2, 3$, with $x_1<x_2<x_3$.
\end{itemize}
 If $\theta=1$ then $f(x)=\frac{L}{4}$ therefore we consider the case $\theta\ne 1$.

 \textbf{3.1. On the set $\mathbb R^{2}_{+}\backslash M_{0}$.}
  Consider operator $W$ on the invariant set $\mathbb R^{2}_{+}\backslash M_{0}$.

To find fixed points of $W$ we have to solve the following system:
\begin{equation}\label{nr}
\begin{array}{l}
      x=L \left(\dfrac{1+x}{1+\theta+x+y}\right)^2\\[2mm]
	  y=L \left(\dfrac{1+y}{1+\theta+x+y}\right)^2.
\end{array}
\end{equation}
Rewrite this as
 $$\sqrt{L}\left(\dfrac{1+x}{1+\theta+x+y}\right)-\sqrt{x}=0$$
$$\sqrt{L}\left(\dfrac{1+y}{1+\theta+x+y}\right)-\sqrt{y}=0.$$

Find $y$ from the first equation:
$$y=\sqrt{L}\left(\dfrac{1}{\sqrt{x}}+\sqrt{x}\right)-(1+\theta+x)=\psi(x).$$

Thus the system of fixed points can be written as:

\begin{equation}\label{aa}
\begin{array}{l}
      y=\psi(x)\\[2mm]
	  x=\psi(y).
\end{array}
\end{equation}
In (\ref{aa})  considen two cases:

{\bf Case}: $x=y$. Then from (\ref{aa}) we get
$$x=\psi(x).$$
That is
$$x=\sqrt{L}\left(\dfrac{1}{\sqrt{x}}+\sqrt{x}\right)-(1+\theta+x).$$

Denoting $u=\sqrt{x}$ we obtain:
$$2u^3-\sqrt{L}u^2+(\theta+1)u-\sqrt{L}=0.$$

Full solution of this equation can be given
by Cardano's formula\footnote{https://en.wikipedia.org/wiki/Cubic$_-$ equation}, they have bulky from:

$u_{1}=\frac{\sqrt{L}}{6}-\frac{6+6\theta-L}{6N}+\frac{1}{6}N,$

$u_{2}=\frac{\sqrt{L}}{6}+\frac{(1+i\sqrt{3})(6+6\theta-L)}{12N}-\frac{1}{12}(1-i\sqrt{3})N,$

$u_{3}=\frac{\sqrt{L}}{6}+\frac{(1-i\sqrt{3})(6+6\theta-L)}{12N}-\frac{1}{12}(1+i\sqrt{3})N,$
 here
$$N=\sqrt[3]{(45-9\theta+L)\sqrt{L}+\sqrt{216+648\theta+648\theta^2+216\theta^3+(1917-1026\theta -27\theta^2 +108L)L}}.$$
 Denote
 $$x^*_1=u_1^2,\ \ x^*_2=u_2^2, \ \ x^*_3=u_3^2.$$

{\bf Case}: $x\neq y$. In this case denote  $u=\sqrt{x}$ and  $v=\sqrt{y}$. Then our system becomes:
$$\left\{\begin{array}{l}
     u^2=\sqrt{L}\left(\dfrac{1}{v}+v\right)-(1+\theta+v^2)  \\
     v^2=\sqrt{L}\left(\dfrac{1}{u}+u\right)-(1+\theta+u^2).
\end{array} \right.$$
 Subtracting we get $(u-v)(\frac{1}{uv}-1)=0.$ The case $u=v$ is considered above. We consider the case $u=\frac{1}{v}$.

Then
$$\frac{1}{u^2}+u^2=\sqrt{L}\left(\dfrac{1}{u}+u\right)-(1+\theta).$$

Denote $\frac{1}{u}+u=\xi$, $\xi\geq2$, then we get:
 $$\xi^2-\sqrt{L}\xi+\theta-1=0.$$
Which has solutions:
 $$\xi_{1,2}=\frac{\sqrt{L}\pm\sqrt{L-4(\theta-1)}}{2}, \ \ \mbox{if} \ \ L\geq4(\theta-1).$$

 Now we find conditions on parameters $\theta,L$ to satisfy $\xi_{1}\geq2, \, \xi_{2}\geq2$. Simple but long analysis show that:

 1) Unique solution to the system (\ref{nr}) exists if parameters are from the set:

 $A_{1,0}=\{(\theta, L):L\leq\frac{(\theta+3)^2}{4},0<\theta\leq5\}\bigcup\{(\theta, L):L<4(\theta-1),\theta>5\}$ (See Figure 1).

 2) Three solutions on set:

 $A_{1,2}=\{(\theta, L):L>\frac{(\theta+3)^2}{4},0<\theta\leq17\}\bigcup\{(\theta, L):\frac{(\theta+3)^2}{4}\leq L<L_{1}, 17<\theta\leq9+8\sqrt{2}\}\bigcup\{(\theta, L):L>L_{2},\theta>17\}\bigcup\{(\theta, L):L=4(\theta-1),\theta>5\}$ (see Figure 2).

 3) Four solutions on set:

 $A_{2,2}=\{(\theta, L):L=L_{2},\theta>17\}\bigcup\{(\theta, L):L=L_{1},17<\theta<9+8\sqrt{2}\}$ (see Figure 3 and Figure 4).

 4) Five solutions on sets:

 $A_{1,4}=\{(\theta, L):4(\theta-1)<L<\frac{(\theta+3)^2}{4},5<\theta\leq9+8\sqrt{2}\}\bigcup\{(\theta, L):4(\theta-1)<L<L_{1},\theta>9+8\sqrt{2}\}$ (see Figure 5), and

 $A_{3,2}=\{(\theta, L):L_{1}<L<L_{2},17<\theta\leq9+8\sqrt{2}\}\bigcup\{(\theta, L):\frac{(\theta+3)^2}{4}<L<L_{2},\theta>9+8\sqrt{2}\}$
 (see Figure 6).

 5) Six solutions on set:

 $A_{2,4}=\{(\theta, L):L=L_{1},\theta>9+8\sqrt{2}\}$ (see Figure 7).

 6) Seven solutions on set:

 $A_{3,4}=\{(\theta, L):L_{1}<L<\frac{(\theta+3)^2}{4},\theta>9+8\sqrt{2}\}$ (see Figure 8).

 \begin{center}

   \begin {tikzpicture} [scale=0.55]
   \draw[->, thick] (0,0) -- (0,10);
   \draw[->, thick] (0,0) -- (11,0);
   \draw[] (1,-0.1) -- (1,0.1);
   \draw[] (2,-0.1) -- (2,0.1);
   \draw[] (3,-0.1) -- (3,0.1);
   \draw[] (4,-0.1) -- (4,0.1);
   \draw[] (5,-0.1) -- (5,0.1);
   \draw[] (6,-0.1) -- (6,0.1);
   \draw[] (7,-0.1) -- (7,0.1);
   \draw[] (8,-0.1) -- (8,0.1);
   \draw[] (9,-0.1) -- (9,0.1);
   \draw[] (10,-0.1) -- (10,0.1);
   \draw[] (-0.1,1) -- (0.1,1);
   \draw[] (-0.1,2) -- (0.1,2);
   \draw[] (-0.1,3) -- (0.1,3);
   \draw[] (-0.1,4) -- (0.1,4);
   \draw[] (-0.1,5) -- (0.1,5);
   \draw[] (-0.1,6) -- (0.1,6);
   \draw[] (-0.1,7) -- (0.1,7);
   \draw[] (-0.1,8) -- (0.1,8);
   \draw[] (-0.1,9) -- (0.1,9);
    \draw[thick, black] (0,0)-- (11,11);
   \draw[thick, red] (0.4,11) .. controls (0.6,0.9).. (0.7,0);
   \draw[thick, green] (11,0.4) .. controls (0.9,0.6).. (0,0.7);

   \node[above] at (11,0){$x$};

   \node[right] at (-0.2,10){$y$};

   \node[left] at (1.5,1.4){$p^{*}_1$};

   \draw[->, thick] (12,0) -- (12,10);
   \draw[->, thick] (12,0) -- (23,0);
   \draw[] (13,-0.1) -- (13,0.1);
   \draw[] (14,-0.1) -- (14,0.1);
   \draw[] (15,-0.1) -- (15,0.1);
   \draw[] (16,-0.1) -- (16,0.1);
   \draw[] (17,-0.1) -- (17,0.1);
   \draw[] (18,-0.1) -- (18,0.1);
   \draw[] (19,-0.1) -- (19,0.1);
   \draw[] (20,-0.1) -- (20,0.1);
   \draw[] (21,-0.1) -- (21,0.1);
   \draw[] (22,-0.1) -- (22,0.1);
   \draw[] (11.9,1) -- (12.1,1);
   \draw[] (11.9,2) -- (12.1,2);
   \draw[] (11.9,3) -- (12.1,3);
   \draw[] (11.9,4) -- (12.1,4);
   \draw[] (11.9,5) -- (12.1,5);
   \draw[] (11.9,6) -- (12.1,6);
   \draw[] (11.9,7) -- (12.1,7);
   \draw[] (11.9,8) -- (12.1,8);
   \draw[] (11.9,9) -- (12.1,9);
    \draw[thick, black] (12,0)-- (23,11);
   \draw[thick, red] (12.5,11).. controls (12.52,0.3) .. (13, 1.5) .. controls (15.3,4) and (15,3) .. (21,0);
   \draw[thick, green] (23 ,0.5).. controls (12.3,0.52) .. (13.5,1) .. controls (16,3.3) and (15,3) .. (12,9);

   \node[above] at (23,0){$x$};

   \node[right] at (11.8,10){$y$};
   \node[right] at (12.7,8){$p_4$};
   \node[right] at (20,0.7){$p_5$};
   \node[left] at (15,3.1){$p^{*}_1$};
\end{tikzpicture}
$\begin{array}{cc}
     \text{Figure 1.}  & \text{Figure 2.} \\
 \text{The case } \theta =20,L=59;
 (20,59)\in A_{1,0} & \text{The case } \theta=6,L=90;
 (6,90)\in A_{1,2} \\
 \end{array} $
\end{center}

\vspace{-0.7 cm}

\begin{center}

   \begin {tikzpicture} [scale=0.55]
   \draw[->, thick] (0,0) -- (0,10);
   \draw[->, thick] (0,0) -- (11,0);
   \draw[] (1,-0.1) -- (1,0.1);
   \draw[] (2,-0.1) -- (2,0.1);
   \draw[] (3,-0.1) -- (3,0.1);
   \draw[] (4,-0.1) -- (4,0.1);
   \draw[] (5,-0.1) -- (5,0.1);
   \draw[] (6,-0.1) -- (6,0.1);
   \draw[] (7,-0.1) -- (7,0.1);
   \draw[] (8,-0.1) -- (8,0.1);
   \draw[] (9,-0.1) -- (9,0.1);
   \draw[] (10,-0.1) -- (10,0.1);
   \draw[] (-0.1,1) -- (0.1,1);
   \draw[] (-0.1,2) -- (0.1,2);
   \draw[] (-0.1,3) -- (0.1,3);
   \draw[] (-0.1,4) -- (0.1,4);
   \draw[] (-0.1,5) -- (0.1,5);
   \draw[] (-0.1,6) -- (0.1,6);
   \draw[] (-0.1,7) -- (0.1,7);
   \draw[] (-0.1,8) -- (0.1,8);
   \draw[] (-0.1,9) -- (0.1,9);
    \draw[thick, black] (0,0)-- (11,11);
   \draw[thick, red] (0.5,11).. controls (0.52,0.3) .. (2.2, 1.8) .. controls (3.3,4) and (3,3) .. (9,0);
   \draw[thick, green] (11 ,0.5).. controls (0.3,0.52) .. (1.8,2.2) .. controls (4,3.3) and (3,3) .. (0,9);

   \node[above] at (11,0){$x$};

   \node[right] at (-0.2,10){$y$};
\node[right] at (0.7,8){$p_4$};
   \node[right] at (8,0.7){$p_5$};
   \node[left] at (1.7,0.7){$p^{*}_1$};
   \node[left] at (2.7,2.8){$p^{*}_2$};

   \draw[->, thick] (12,0) -- (12,10);
   \draw[->, thick] (12,0) -- (23,0);
   \draw[] (13,-0.1) -- (13,0.1);
   \draw[] (14,-0.1) -- (14,0.1);
   \draw[] (15,-0.1) -- (15,0.1);
   \draw[] (16,-0.1) -- (16,0.1);
   \draw[] (17,-0.1) -- (17,0.1);
   \draw[] (18,-0.1) -- (18,0.1);
   \draw[] (19,-0.1) -- (19,0.1);
   \draw[] (20,-0.1) -- (20,0.1);
   \draw[] (21,-0.1) -- (21,0.1);
   \draw[] (22,-0.1) -- (22,0.1);
    \draw[] (11.9,1) -- (12.1,1);
   \draw[] (11.9,2) -- (12.1,2);
   \draw[] (11.9,3) -- (12.1,3);
   \draw[] (11.9,4) -- (12.1,4);
   \draw[] (11.9,5) -- (12.1,5);
   \draw[] (11.9,6) -- (12.1,6);
   \draw[] (11.9,7) -- (12.1,7);
   \draw[] (11.9,8) -- (12.1,8);
   \draw[] (11.9,9) -- (12.1,9);
    \draw[thick, black] (12,0)-- (23,11);
   \draw[thick, red] (12.5,11).. controls (12.52,0.7) .. (14, 2) .. controls (15.2,4.7) and (15,3) .. (21,0);
   \draw[thick, green] (23 ,0.5).. controls (12.7,0.52) .. (14,2) .. controls (16.7,3.2) and (15,3) .. (12,9);

   \node[above] at (23,0){$x$};

   \node[right] at (11.8,10){$y$};
    \node[right] at (12.7,8){$p_4$};
   \node[right] at (20,0.7){$p_5$};
   \node[left] at (14.3,2.4){$p^{*}_1$};
   \node[right] at (15.5,3.4){$p^{*}_2$};

\end{tikzpicture}

$\begin{array}{cc}
     \text{Figure 3.}  & \text{Figure 4.} \\
 \text{The case } \theta =19,L=125;
 (19,125)\in A_{2,2} & \text{The case } \theta=19,L=128;
 (19,128)\in A_{2,2} \\
 \end{array} $
\end{center}
\vspace{-0.4 cm}
\begin{center}

   \begin {tikzpicture} [scale=0.55]
   \draw[->, thick] (0,0) -- (0,10);
   \draw[->, thick] (0,0) -- (11,0);
   \draw[] (1,-0.1) -- (1,0.1);
   \draw[] (2,-0.1) -- (2,0.1);
   \draw[] (3,-0.1) -- (3,0.1);
   \draw[] (4,-0.1) -- (4,0.1);
   \draw[] (5,-0.1) -- (5,0.1);
   \draw[] (6,-0.1) -- (6,0.1);
   \draw[] (7,-0.1) -- (7,0.1);
   \draw[] (8,-0.1) -- (8,0.1);
   \draw[] (9,-0.1) -- (9,0.1);
   \draw[] (10,-0.1) -- (10,0.1);
   \draw[] (-0.1,1) -- (0.1,1);
   \draw[] (-0.1,2) -- (0.1,2);
   \draw[] (-0.1,3) -- (0.1,3);
   \draw[] (-0.1,4) -- (0.1,4);
   \draw[] (-0.1,5) -- (0.1,5);
   \draw[] (-0.1,6) -- (0.1,6);
   \draw[] (-0.1,7) -- (0.1,7);
   \draw[] (-0.1,8) -- (0.1,8);
   \draw[] (-0.1,9) -- (0.1,9);
    \draw[thick, black] (0,0)-- (11,11);
   \draw[thick, red] (0.5,11).. controls (0.52,0.3) .. (1, -0.8) .. controls (3.3,4) and (3,3) .. (9,0);
   \draw[thick, green] (11 ,0.5).. controls (0.3,0.52) .. (-0.8,1) .. controls (4,3.3) and (3,3) .. (0,9);

   \node[above] at (11,0){$x$};

   \node[right] at (-0.2,10){$y$};
\node[right] at (0.7,8){$p_3$};
   \node[right] at (8,0.7){$p_4$};
   \node[left] at (1.7,0.8){$p^{*}_1$};
   \node[left] at (1.2,2.1){$p_1$};
   \node[left] at (2.8,0.8){$p_2$};

   \draw[->, thick] (12,0) -- (12,10);
   \draw[->, thick] (12,0) -- (23,0);
   \draw[] (13,-0.1) -- (13,0.1);
   \draw[] (14,-0.1) -- (14,0.1);
   \draw[] (15,-0.1) -- (15,0.1);
   \draw[] (16,-0.1) -- (16,0.1);
   \draw[] (17,-0.1) -- (17,0.1);
   \draw[] (18,-0.1) -- (18,0.1);
   \draw[] (19,-0.1) -- (19,0.1);
   \draw[] (20,-0.1) -- (20,0.1);
   \draw[] (21,-0.1) -- (21,0.1);
   \draw[] (22,-0.1) -- (22,0.1);
     \draw[] (11.9,1) -- (12.1,1);
   \draw[] (11.9,2) -- (12.1,2);
   \draw[] (11.9,3) -- (12.1,3);
   \draw[] (11.9,4) -- (12.1,4);
   \draw[] (11.9,5) -- (12.1,5);
   \draw[] (11.9,6) -- (12.1,6);
   \draw[] (11.9,7) -- (12.1,7);
   \draw[] (11.9,8) -- (12.1,8);
   \draw[] (11.9,9) -- (12.1,9);
    \draw[thick, black] (12,0)-- (23,11);
   \draw[thick, red] (12.5,11).. controls (12.52,0.3) .. (14.2, 1.8) .. controls (15.2,4.7) and (15,3) .. (21,0);
   \draw[thick, green] (23 ,0.5).. controls (12.3,0.52) .. (13.8,2.2) .. controls (16.7,3.2) and (15,3) .. (12,9);

   \node[above] at (23,0){$x$};

   \node[right] at (11.8,10){$y$};
    \node[right] at (12.7,8){$p_4$};
   \node[right] at (20,0.7){$p_5$};
   \node[left] at (13.7,0.6){$p^{*}_1$};
   \node[left] at (14.7,2.8){$p^{*}_2$};
   \node[right] at (15.5,3.4){$p^{*}_3$};

\end{tikzpicture}

$\begin{array}{cc}
     \text{Figure 5.}  & \text{Figure 6.} \\
 \text{The case } \theta =19,L=98;
 (19,98)\in A_{1,4} & \text{The case } \theta=19,L=126;
 (19,126)\in A_{3,2} \\
 \end{array} $
\end{center}

\begin{center}

   \begin {tikzpicture} [scale=0.55]
   \draw[->, thick] (0,0) -- (0,10);
   \draw[->, thick] (0,0) -- (11,0);
   \draw[] (1,-0.1) -- (1,0.1);
   \draw[] (2,-0.1) -- (2,0.1);
   \draw[] (3,-0.1) -- (3,0.1);
   \draw[] (4,-0.1) -- (4,0.1);
   \draw[] (5,-0.1) -- (5,0.1);
   \draw[] (6,-0.1) -- (6,0.1);
   \draw[] (7,-0.1) -- (7,0.1);
   \draw[] (8,-0.1) -- (8,0.1);
   \draw[] (9,-0.1) -- (9,0.1);
   \draw[] (10,-0.1) -- (10,0.1);
   \draw[] (-0.1,1) -- (0.1,1);
   \draw[] (-0.1,2) -- (0.1,2);
   \draw[] (-0.1,3) -- (0.1,3);
   \draw[] (-0.1,4) -- (0.1,4);
   \draw[] (-0.1,5) -- (0.1,5);
   \draw[] (-0.1,6) -- (0.1,6);
   \draw[] (-0.1,7) -- (0.1,7);
   \draw[] (-0.1,8) -- (0.1,8);
   \draw[] (-0.1,9) -- (0.1,9);
    \draw[thick, black] (0,0)-- (11,11);
  \draw[thick, red] (0.5,11).. controls (0.52,0.3) .. (1, 0.2) .. controls (3,4) and (3,3) .. (9,0);
   \draw[thick, green] (11 ,0.5).. controls (0.3,0.52) .. (0.2,1) .. controls (4,3) and (3,3) .. (0,9);

   \node[above] at (11,0){$x$};

   \node[right] at (-0.2,10){$y$};
\node[right] at (0.7,8){$p_3$};
   \node[right] at (8,0.7){$p_4$};
   \node[left] at (1.7,0.8){$p^{*}_1$};
   \node[left] at (1.4,2.1){$p_1$};
   \node[left] at (2.8,0.8){$p_2$};
   \node[left] at (2.5,2.5){$p^{*}_2$};

   \draw[->, thick] (12,0) -- (12,10);
   \draw[->, thick] (12,0) -- (23,0);
   \draw[] (13,-0.1) -- (13,0.1);
   \draw[] (14,-0.1) -- (14,0.1);
   \draw[] (15,-0.1) -- (15,0.1);
   \draw[] (16,-0.1) -- (16,0.1);
   \draw[] (17,-0.1) -- (17,0.1);
   \draw[] (18,-0.1) -- (18,0.1);
   \draw[] (19,-0.1) -- (19,0.1);
   \draw[] (20,-0.1) -- (20,0.1);
   \draw[] (21,-0.1) -- (21,0.1);
   \draw[] (22,-0.1) -- (22,0.1);
     \draw[] (11.9,1) -- (12.1,1);
   \draw[] (11.9,2) -- (12.1,2);
   \draw[] (11.9,3) -- (12.1,3);
   \draw[] (11.9,4) -- (12.1,4);
   \draw[] (11.9,5) -- (12.1,5);
   \draw[] (11.9,6) -- (12.1,6);
   \draw[] (11.9,7) -- (12.1,7);
   \draw[] (11.9,8) -- (12.1,8);
   \draw[] (11.9,9) -- (12.1,9);
    \draw[thick, black] (12,0)-- (23,11);
   \draw[thick, red] (12.5,11).. controls (12.52,0.3) .. (13, 0.2) .. controls (14.4,4) and (15,3) .. (21,0);
   \draw[thick, green] (23 ,0.5).. controls (12.3,0.52) .. (12.2,1) .. controls (16,2.4) and (15,3) .. (12,9);

   \node[above] at (23,0){$x$};

   \node[right] at (11.8,10){$y$};
    \node[right] at (12.7,8){$p_3$};
   \node[right] at (20,0.7){$p_4$};
   \node[left] at (13.3,0.3){$p^{*}_1$};
   \node[left] at (13.1,1.9){$p_1$};
   \node[left] at (14.4,0.8){$p_2$};
   \node[left] at (13.9,1.8){$p^{*}_2$};
   \node[left] at (14.6,2.7){$p^{*}_3$};

\end{tikzpicture}

$\begin{array}{cc}
     \text{Figure 7.}  & \text{Figure 8.} \\
 \text{The case } \theta =22,5;L=154;
 (22,5;154)\in A_{2,4} & \text{The case } \theta=22,5;L=157;
 (22,5;157)\in A_{3,4} \\
 \end{array} $

\end{center}

 Denoting $A_{i,j}$ above means that the number of fixed points on $M_{0}$ is $i$ and number of fixed points outside $M_{0}$ is $j$. Thus we proved the following:

\textbf{Lemma 4.} The following assertions hold
\begin{itemize}
	\item[1)] If $(\theta, L)\in A_{1,0}$ then

 $Fix (W)=\{p^*_{1}\}$ .
	
	\item[2)] If $(\theta, L)\in A_{1,2}$ then

 $Fix(W)=\{p^*_{1}, p_{1}, p_{2}\}$ .

    \item[3)] If $(\theta, L)\in A_{2,2}$ then

    $Fix(W)=\{p^*_{1}, p^*_{2}, p_{1}, p_{2}\}$ .
		
    \item[4)] If $(\theta, L)\in A_{1,4}$ then

     $Fix(W)=\{p^*_{1}, p_{1}, p_{2}, p_{3}, p_{4}\}$ .

    \item[5)] If $(\theta, L)\in A_{3,2}$ then

     $Fix(W)=\{p^*_{1}, p^*_{2}, p^*_{3}, p_{1}, p_{2}\}$ .

     \item[6)] If $(\theta, L)\in A_{2,4}$ then

     $Fix(W)=\{p^*_{1}, p^*_{2}, p_{1}, p_{2}, p_{3}, p_{4}\}$ .

    \item[7)] If $(\theta, L)\in A_{3,4}$ then

     $Fix(F)=\{p^*_{1}, p^*_{2}, p^*_{3}, p_{1}, p_{2}, p_{3}, p_{4}\}$,
 \end{itemize}
	where
	$$p^*_i=(x^*_i, x^*_i), \, i=1,2,3.$$
		$$p_{1}=(x_{1},x_{2}), p_{2}=(x_{2},x_{1}), p_{3}=(x_{3},x_{4}), p_{4}=(x_{4},x_{3}),$$
with	 $$x_{1}=\frac{\sqrt{L}+\sqrt{L-4(\theta-1)}+\sqrt{2L-4\theta+2\sqrt{L^2-4(\theta-1)L}}}{2},$$
	$$x_{2}=\frac{\sqrt{L}+\sqrt{L-4(\theta-1)}-\sqrt{2L-4\theta+2\sqrt{L^2-4(\theta-1)L}}}{2},$$
	$$x_{3}=\frac{\sqrt{L}-\sqrt{L-4(\theta-1)}+\sqrt{2L-4\theta-2\sqrt{L^2-4(\theta-1)L}}}{2},$$
	$$x_{4}=\frac{\sqrt{L}-\sqrt{L-4(\theta-1)}-\sqrt{2L-4\theta-2\sqrt{L^2-4(\theta-1)L}}}{2}.$$

By Lemma 2 and Lemma 4 we get the following main result.

\textbf{Theorem 1.} Fixed points of operator (\ref{e4o}) for $L_{1}=L_{2}=L$ are as follows:

$$Fix(F)=\left\{\begin{array}{ll}
		\{P_1\}, \ \ \mbox{if} \ \ (\theta,L)\in A_{1,0},\\[2mm]
       \{P_1, P_4, P_5\} \ \ \mbox{if} \ \ (\theta,L)\in A_{1,2},\\[2mm]
       \{P_1, P_2, P_4, P_5\}  \ \ \mbox{if} \ \ (\theta,L)\in A_{2,2},\\[2mm]
        \{P_1, P_4, P_5, P_6, P_7\}  \ \ \mbox{if} \ \ (\theta,L)\in A_{1,4},\\[2mm]
       \{P_1, P_2, P_3, P_4, P_5\}  \ \ \mbox{if} \ \ (\theta,L)\in A_{3,2},\\[2mm]
       \{P_1, P_2, P_4, P_5, P_6, P_7\} \ \ \mbox{if} \ \ (\theta,L)\in A_{2,4},\\[2mm]
		 \{P_1, P_2, P_3, P_4, P_5, P_6, P_7\} \ \ \mbox{if} \ \ (\theta,L)\in A_{3,4};
		\end{array}\right.$$
where
$$\begin{array}{ll}

P_1=(\frac{x^{*}_1}{L}\lambda_1,\frac{x^{*}_1}{L}\lambda_2,\frac{x^{*}_1}{L}\lambda_3,...)\\[3mm]

P_2=(\frac{x^{*}_2}{L}\lambda_1,\frac{x^{*}_2}{L}\lambda_2,\frac{x^{*}_2}{L}\lambda_3,...)\\[3mm]

P_3=(\frac{x^{*}_3}{L}\lambda_1,\frac{x^{*}_3}{L}\lambda_2,\frac{x^{*}_3}{L}\lambda_3,...)\\[3mm]

P_4=(\frac{x_1}{L}\lambda_1,\frac{x_2}{L}\lambda_2,\frac{x_1}{L}\lambda_3,...)\\[3mm]

P_5=(\frac{x_2}{L}\lambda_1,\frac{x_1}{L}\lambda_2,\frac{x_2}{L}\lambda_3,...)\\[3mm]

P_6=(\frac{x_3}{L}\lambda_1,\frac{x_4}{L}\lambda_2,\frac{x_3}{L}\lambda_3,...)\\[3mm]

P_7=(\frac{x_4}{L}\lambda_1,\frac{x_3}{L}\lambda_2,\frac{x_4}{L}\lambda_3,...)

\end{array}$$

\hfill $\Box$\\

\newpage

\begin{center}
	{\sc Cheksiz o'lchamli operatorning qo'zg'almas nuqtalarining tavsifi }\\
	{\bf Olimov Umrbek}\\
\end{center}
\medskip

{Biz sanoqli $\mathbb{N}$ spin qiymatlari to'plamiga ega bo'lgan Hard-Core (HC) modeli bilan bog'liq cheksiz o'lchamli nochiziqli operatorni qaraymiz. Ravshanki, cheksiz o'lchamli operatorning qo'zg'almas nuqtalarini topish umuman olganda imkonsizdir. Lekin biz cheksiz o'lchamli operatorni ikki o'lchamli operatorga keltirish usuli yordamida cheksiz o'lchamli operatorning qo'zg'almas nuqtalarini to'la analiz qildik. Parametrlar to'plamini $A_{i,j}$ qism to'plamlarga ajratdik, bu yerda $i$ indeks $y=x$ to'g'ri chiziqdagi qo'zg'almas nuqtalar sonini hamda $j$ indeks $y=x$ dan tashqaridagi qo'zg'almas nuqtalar sonini ifodalaydi. Qo'zg'almas nuqtalar soni ko'pi bilan yettitagacha bo'lishi mumkin va har bir qo'zg'almas nuqtaning aniq ko'rinishi topilgan.}\\

\medskip
{{\bf Kalit so`zlar:} Qo'zg'almas nuqta,invariant to'plam,cheksiz o'lchamli operator.}\\

\bigskip

\begin{center}
{\sc Описание неподвижных точек бесконечномерного оператора }\\
{\bf Олимов Умрбек}\\
\end{center}
\medskip

{Мы рассматриваем бесконечномерный нелинейный оператор, связанный с моделью Hard-Core (HC) со счетным множеством $\mathbb{N}$ значений спина. Известно, что найти неподвижные точки бесконечномерного оператора вообще говаря невозможно. Но мы полностью проанализировали неподвижные точки бесконечномерного оператора, применив технику сведения бесконечномерного оператора к двумерному оператору. Множество параметров разбито на подмножества $A_{i,j},$, где индекс $i$ означает количество неподвижных точек на прямой $y=x$, $j$ означает количество неподвижных точек вне $ y=x.$ Число неподвижных точек может достигать семи и найден явный вид каждой неподвижной точки.

\medskip
{{\bf Ключевые слова:} Неподвижная точка, инвариантное множество, бесконечномерный оператор.}\\

\bigskip

\textbf{Received: **/**/202*}\\

\textbf{Accepted: **/**/202*}\\

\bigskip
\textcolor{red}{\Large Cite this article}

Olimov U.~ {Description of fixed points of an infinite dimensional operator}  \textit{Bull.~Inst.~Math.},~ 2023, Vol.6, No 1, pp. \pageref{firstpage}-\pageref{lastpage}

\label{lastpage}

\end{document}